\documentclass[12 pt]{amsart}

\usepackage{amsrefs}

\theoremstyle{plain}
\newtheorem{theorem}{Theorem}
\newtheorem{lemma}[theorem]{Lemma}
\newtheorem{proposition}[theorem]{Proposition}

\theoremstyle{definition}
\newtheorem{example}{Example}

\DeclareMathOperator{\Aut}{Aut}
\DeclareMathOperator{\End}{End}
\DeclareMathOperator{\Gal}{Gal}
\DeclareMathOperator{\Tr}{Tr}

\newcommand{\bG}{\mathbb G}
\newcommand{\bZ}{\mathbb Z}
\newcommand{\bQ}{\mathbb Q}
\newcommand{\bR}{\mathbb R}
\newcommand{\bC}{\mathbb C}

\newcommand{\sA}{\mathcal A}

\newcommand{\sG}{\mathcal G}
\newcommand{\sY}{\mathcal Y}

\begin{document}

\title{Hodge structures of CM-type}

\author{Salman Abdulali}
\address{Department of Mathematics, East Carolina University, 
Greenville, NC 27858, USA}
\email{abdulalis@mail.ecu.edu}
\thanks{Research supported in part by NSA grant number H98230-04-1-0102, and in part by a Research/Creative Activity Grant from East Carolina University.}

\subjclass[2000]{Primary 14C30, 14K22}

\begin{abstract}
We show that any effective Hodge structure of CM-type occurs (without having to take a Tate twist) in the cohomology of some CM abelian variety over $\bC$.
As a consequence we get a simple proof of the theorem (due to Hazama) that the usual Hodge conjecture for the class of all CM abelian varieties implies the general Hodge conjecture for the same class.
\end{abstract}

\maketitle

\section{Introduction}
Recall that a Hodge structure $V_{\bC} = \bigoplus_{p+q=n} V^{p,q}$ is said to be effective if $V^{p,q} = 0$ unless $p,q \geq 0$, and, it is said to be geometric if it is isomorphic to a Hodge substructure of the cohomology of a smooth, projective variety over $\bC$.
For $m \in \bZ$, the Tate twist $V(m)$ is defined by $V(m)^{p,q} = V^{p+m,q+m}$.

A geometric Hodge structure must be effective, though not conversely.
The general Hodge conjecture as formulated by Grothendieck \cite{Grothendieck} implies that any effective Tate twist of a geometric Hodge structure is again geometric.

A Hodge structure $V$ is said to be of \emph{CM-type} if it is polarizable, and, its Hodge group $G$ is abelian.
In this paper we show (Theorem \ref{main}) that any effective Hodge structure of CM-type is geometric.
Moreover, such a Hodge structure is isomorphic (without having to take a Tate twist) to a Hodge structure appearing in the cohomology of a CM abelian variety.
That this is true up to a Tate twist is due to Serre (see \cite{MilneShih}*{(1.7), p.~234}).

In a series of papers \cites{Abdulali1997, Abdulali2000, Abdulali2001, Abdulali2002, Abdulali2004}, we have shown that for certain abelian varieties $A$ over $\bC$, the general Hodge conjecture for $A$ is implied by the usual Hodge conjecture for a certain class of abelian varieties.
Chad Schoen and Fumio Hazama have kindly pointed out to me an error in the proof of Prop.~4.4.1 of \cite{Abdulali1997}.
Consequently, this proposition, Theorem 6.1 of \cite{Abdulali1997}, and, Theorem 3.1 of \cite{Abdulali2001}, are false as stated.
Correcting these errors, we prove in an elementary manner that the usual Hodge conjecture for all CM abelian varieties implies the general Hodge conjecture for the same class.
This result has been independently obtained by Hazama \cites{Hazama2002, Hazama2003} using different methods.
Finally, we give several examples of CM abelian varieties for which the general Hodge conjecture can be proved by these methods.

\section{Hodge structures}
By a \emph{rational Hodge structure} of \emph{weight} $n$ we mean a finite dimensional vector space $V$ over $\bQ$, and a decomposition
$V_{\bC} = \bigoplus_{p+q=n} V^{p,q}$
such that $\overline{V^{p,q}} = V^{q,p}$ for all $p,q$.
It is \emph{effective} if $V^{p,q} = 0$ unless $p,q \geq 0$.
A Hodge structure $V$ determines a morphism
$$\lambda_V : S^1 \to GL(V_{\bR}),$$
where $S^1$ is the unit circle in the complex plane, such that $\lambda_V(e^{i\theta})$ acts on $V^{p,q}$ as multiplication by $e^{(p-q)i\theta}$.
The element $C = \lambda_V(i)$ is called the \emph{Weil operator}.
A \emph{polarization} of $V$ is a morphism of Hodge structures $\psi : V \otimes V \to \bQ(-n)$ such that $\psi (x,Cy)$ is symmetric and positive definite on $V_{\bR}$.

The \emph{Hodge group} (or Special Mumford-Tate group) of $V$ is defined to be the smallest algebraic subgroup $G$ of $GL(V)$ such that $G(\bR)$ contains the image of $\lambda_V$.
It is a reductive algebraic group over $\bQ$ if $V$ is polarizable \cite{900}*{Prop.~3.6, p.~44}.

Let $A$ be an abelian variety over $\bC$.
The Hodge group of $A$ is the Hodge group of $H^1(A,\bQ)$.
$A$ is of CM-type if and only if its Hodge group is abelian (Mumford \cite{Mumford1969}*{p.~347}).
A rational Hodge structure $V$ is said to be of \emph{CM-type} if it is polarizable, and, its Hodge group $G$ is abelian.
It is said to be \emph{trivial} if $G$ is trivial, and \emph{nontrivial} otherwise.

We shall now collect together some facts about CM Hodge structures.
For the convenience of the reader we include proofs of some \lq\lq{well-known}\rq\rq\ facts.
The main references are Deligne \cite{900} and Schappacher \cite{Schappacher}.

Let $V$ be a Hodge structure with Hodge group $G$.
The endomorphism algebra of $V$ is given by $D := \End_G V$.
If $E$ is a subfield of $D$, then, $V$ may be considered as a vector space over $E$.
We say that $V$ is of \emph{type} $E$ if $V$ is $1$-dimensional over $E$.

Assume that $V$ is an irreducible Hodge structure of CM-type.
Then $D$ is a division algebra.
Let $E$ be the center of $D$, and let $T_E$ be the torus $R_{E/\bQ} \bG_{m,E}$.
Since $G$ is abelian, it is contained in $T_E$.
Consider $V$ as a vector space over $E$.
Any $E$-subspace of $V$ is also a $G$-submodule of $V$.
Since $V$ is irreducible, it must be $1$-dimensional over $E$.
Therefore $D = \End_G V = E$.
Further, from the irreducibility of $V$, we see that $G(\bQ)$ generates $E$ as an algebra over $\bQ$.

We have $E_{\bC} := E \otimes_{\bQ} \bC = \oplus_{\sigma\in S} E^{\sigma}$, where $E^{\sigma} = E \otimes_{E,\sigma} \bC$, and, $S$ is the set of embeddings of $E$ into $\bC$.
We consider $E$ embedded in $E_{\bC}$ via $e \mapsto (\sigma(e))$.
Then the trace map $\Tr_{E_{\bC}/\bC}: E_{\bC} \to \bC$ given by
$$\bigoplus_{\sigma\in S} E^{\sigma} \ni (a_{\sigma}) \mapsto \sum_{\sigma \in S} a_{\sigma}$$
is compatible with $\Tr_{E/{\bQ}}: E \to \bQ$.

\begin{lemma}
If $V$ is a nontrivial irreducible Hodge structure of CM-type with Hodge group $G$, then, $E = \End_G V$ is a CM-field.
\end{lemma}

\begin{proof}
Let $\psi: V \times V \to \bQ$ be a polarization of $V$.
Since $G$ is a torus contained in $\Aut \psi$, which is either an orthogonal group or a symplectic group, we have that $G(\bR)$ is compact.
It follows that $E$ is totally imaginary, for a real embedding of $E$ would induce a nontrivial character of $G(\bR)$ with noncompact image.

The polarization $\psi$ determines an involution $e \mapsto e^{\star}$ of $E$ by the rule
$$\psi(x,ey) = \psi(e^{\star}x,y).$$
Let $F$ be the fixed field of this involution.
Then $[E:F]$ is either $1$ or $2$, depending on whether the involution is trivial or not.
For $g \in G$ we have $\psi(x,y) = \psi(gx,gy) = \psi(g^{\star}gx,y)$, so $g^{\star}g = 1$.
Thus $\star$ is nontrivial, and, $[E:F] = 2$.
We extend $\star$ to an involution of $E \otimes_{\bQ} \bR$; the set of elements fixed by the extended involution is $F \otimes_{\bQ} \bR$.

Let $n$ be the weight of the Hodge structure.
Since the Weil operator $C$ is an element of $G(\bR)$ such that $C^2 = (-1)^n$, we have $C^{\star} = C^{-1} = (-1)^nC$.

Suppose $n$ is even.
There exists a unique form $T: V \times V \to E$ such that $T(y,x) = T(x,y)^{\star}$,
$T(ax,by) = ab^{\star}T(x,y)$, and, $\psi(x,y) = \Tr_{E/\bQ} T(x,y)$ for $a,b \in E$, $x,y \in V$.
Identify $V$ with $E$, and let $\beta = T(1,1)$.
Then $\psi(x,Cx) = \Tr_{E \otimes_{\bQ} \bR/\bR} T(x,Cx) = \Tr_{E \otimes_{\bQ} \bR/\bR} C\beta xx^{\star} > 0$ for $x \neq 0$.
Note that $f := C\beta \in F \otimes_{\bQ} \bR$.

Suppose next that $n$ is odd.
Let $\alpha \in F$ be such that $E = F(\sqrt{\alpha})$.
Then $\sqrt{\alpha}^{\,\star} = -\sqrt{\alpha}$.
There exists a unique form $T: V \times V \to E$ such that $T(y,x) = -T(x,y)^{\star}$,
$T(ax,by) = ab^{\star}T(x,y)$, and, $\psi(x,y) = \Tr_{E/\bQ} \sqrt{\alpha} T(x,y)$ for $a,b \in E$, $x,y \in V$.
Identify $V$ with $E$, and let $\beta = T(1,1)$.
Then
$$\psi(x,Cx) = \Tr_{E \otimes_{\bQ} \bR/\bR} \sqrt{\alpha} T(x,Cx) = \Tr_{E \otimes_{\bQ} \bR/\bR} C^{\star} \sqrt{\alpha} \beta xx^{\star} > 0$$ for $x \neq 0$.
Note that $f := C^{\star} \sqrt{\alpha} \beta \in F \otimes_{\bQ} \bR$.

Thus in all cases we have $\Tr_{E\otimes_{\bQ} \bR/\bR} fxx^{\star} > 0$ for all $x \in (E \otimes_{\bQ} \bR)^{\times}$, and a fixed $f \in F \otimes_{\bQ} \bR$.
In particular, for all $x \in (F \otimes_{\bQ} \bR)^{\times}$ we have $\Tr_{F\otimes_{\bQ} \bR/\bR} fx^2 > 0$.
Write $f = (f_{\sigma}) \in \oplus_{\sigma \in S} E^{\sigma}$.

Suppose $\tau: F \to \bC$ were a nonreal embedding.
Write $$f_{\sigma} = r_{\sigma} e^{i\theta_{\sigma}} \otimes s_{\sigma} \in \bC \otimes_{\bQ} \bR,$$
with all $r_{\sigma}, s_{\sigma} > 0$.
By the Artin-Whaples approximation theorem there exists $x \in F$ such that $\sigma(x)$ is close to $0$ for all embeddings of $F$ except $\tau$ and $\overline{\tau}$, while $\tau(x)$ has large absolute value, and argument close to $\frac{1}{2} (\pi - \theta_{\sigma})$.
Then $\Tr_{F\otimes_{\bQ} \bR/\bR} fx^2 < 0$.
This contradiction shows that $F$ is totally real.
\end{proof}

\begin{lemma}
\label{typelemma}
Let $E$ be a CM-field, and $S$ the set of embeddings of $E$ into $\bC$.
Let $V$ be a $1$-dimensional vector space over $E$.
There is a bijection $V_{\phi} \leftrightarrow \phi$ between Hodge structures on $V$ of CM-type of weight $n$ with endomorphisms by $E$, and, functions $\phi: S \to \bZ$ satisfying
\begin{equation}
\label{type}
\phi(\sigma) + \phi(\bar{\sigma}) = n.
\end{equation}
\end{lemma}

\begin{proof}
Let $V$ be a Hodge structure which is $1$-dimensional as a vector space over $E$.
Let $G$ be the Hodge group of $V$.
Then,
$$V \otimes_{\bQ} \bC = E \otimes_{\bQ} \bC = \bigoplus_{\sigma \in S} V^{\sigma},$$
where, $V^{\sigma} := E \otimes_{E,\sigma} \bC = E^{\sigma}$.
Since $G \subset T_E$, each $V^{\sigma}$ is a $G_{\bC}$-submodule of $V_\bC$ of complex dimension $1$.
Hence $V^{\sigma} \subset V^{p,q}$ for some $p,q$, with $p+q=n$.
Write $\phi(\sigma) = p$.
Since $\overline{V^{\sigma}} = V^{\bar{\sigma}}$, $\phi$ satisfies \eqref{type}.

Conversely, given a function $\phi: S \to \bZ$ satisfying \eqref{type}, we can define a unique Hodge structure $V_{\phi}$ on $V$ such that $V^{\sigma}$ has Hodge type $(\phi(\sigma), \phi(\bar{\sigma}))$.
Let $H$ be the Hodge group of $V_{\phi}$.
Since each $V^{\sigma}$ is an $E$-submodule of $V_{\bC}$, we see that $\End_H V_{\phi} \supset E$.
Then $E$ equals its own centralizer in $\End_{\bQ} V$, so $H \subset T_E$, and $V_{\phi}$ is of CM-type.
\end{proof}

\section{Abelian varieties}

Let $X$ be a smooth projective variety over $\bC$, and $\sY$ a family of smooth projective varieties over $\bC$.
We say that $X$ is \emph{dominated} by $\sY$ if, for each irreducible Hodge substructure $V$ of the cohomology of $X$, there exists $Y \in \sY$, a nonnegative integer $n$, and a Hodge substructure $W \subset H^n(Y, \bQ)$ such that $W$ is isomorphic to a Tate twist of $V$, and, $W_{\bC}$ contains an element of Hodge type $(n,0)$.

If $X$ is dominated by $\sY$, then, the usual Hodge conjecture for all $X \times Y$ for all $Y \in \sY$ implies the general Hodge conjecture for $X$.
This observation is due to Grothendieck; see \cite{Abdulali1997}*{Prop.~2.1, p.~243} for a proof.

\begin{theorem}
\label{main}
Let $V$ be an effective Hodge structure with complex multiplication by $E$.
Then $V$ is isomorphic to a Hodge substructure of the cohomology of a product of abelian varieties with CM by $E$.
\end{theorem}

\begin{proof}
We shall prove this by induction on the weight $n$ of the Hodge structure $V$, the case $n=1$ being true by the classical theory of complex multiplication.
Let $\phi: S \to \bZ$ be the function associated to $V$ by Lemma \ref{typelemma}, so that $V = V_{\phi}$.
Choose a CM-type $T$ of $E$ such that $\phi(\sigma) \geq \phi(\bar{\sigma})$ for all $\sigma \in T$.
Let $\chi$ be the characteristic function of $T$.
Then $V_{\chi} = H^1(A, \bQ)$ for an abelian variety $A$ with CM by $E$.
Let $\psi = \phi - \chi$.
Then $V_{\psi}$ is an effective Hodge structure of weight $n-1$, so by induction, it is contained in $H^r(B, \bQ)$, with $B$ a product of abelian varieties with CM by $E$.
Inside the K\"unneth component $H^1(A,\bC) \otimes H^r(B, \bC)$ of $H^{r+1}(A \times B, \bC)$, we have the space
$$\bigoplus_{\sigma \in S} V_{\chi}^{\sigma} \otimes V_{\psi}^{\sigma} = \bigoplus_{\sigma \in S} V_{\chi + \psi}^{\sigma} = \bigoplus_{\sigma \in S} V^{\sigma} = V_{\bC},$$
so, $H^{r+1}(A \times B, \bQ)$ contains a Hodge substructure isomorphic to $V$.
\end{proof}

\begin{theorem}
Any abelian variety of CM-type is dominated by the class of all CM abelian varieties.
Let $E$ be a CM-field, and $A$ an abelian variety with CM by $E$.
Then any power of $A$ is dominated by the set of products of all CM abelian varieties with CM-field $E$.
\end{theorem}

\begin{proof}
Let $A$ be a CM abelian variety, and, $V$ an irreducible Hodge structure in $H^n(A, \bQ)$.
Then $V$ is itself of CM-type since its Hodge group is a quotient of the Hodge group of $A$ (\cite{Abdulali2002}*{Lemma~2.1.1, p.~917}).
Let $r$ be the largest integer such that $V(r)$ is effective.
By the previous theorem, $V(r)$ occurs in the cohomology of a product of CM abelian varieties.
This shows that $A$ is dominated by the class of all CM abelian varieties.

Suppose now that $A$ has CM by a CM-field $E$.
We assume without loss of generality that $A$ is simple, and identify $U := H^1(A,\bQ)$ with $E$.
Let $G$ be the Hodge group of $A$.
Let $V$ be an irreducible Hodge structure in the cohomology of a power of $A$.
Then $V$ is isomorphic to a Hodge substructure of $U^{\otimes n}$ for a positive integer $n$.
We note that the Hodge group of $U^{\otimes n}$ is also $G$.
The identification of $U$ with $E$ extends to an identification of $U^{\otimes n}$ with $E^{\otimes n}$.
Since $G \subset E^{\otimes n}$, any ideal of $E^{\otimes n}$ is also a $G$-submodule, and hence a Hodge structure.
Hence $V$ is contained in a simple ideal $W$ of $E^{\otimes n}$.
But $E^{\otimes n}$ is isomorphic to a direct sum of copies of $E$.
Hence $W$ is a Hodge structure of type $E$.
It follows that $F$, the CM-field of $V$, is contained in $E$.
Let $r$ be the largest integer such that $V(r)$ is effective.
By the previous theorem, $V(r)$ occurs in the cohomology of a product of abelian varieties with CM by $F$.
Since any CM-type for $F$ can be lifted to a CM-type for $E$, $V(r)$ also occurs in the cohomology of a product of abelian varieties with CM by $E$, and $A$ is dominated by such products.
\end{proof}

\begin{proposition}
\label{product}
Let $E_1, \dots, E_m$ be CM-fields whose Galois closures are linearly disjoint over $\bQ$.
For each $i=1,\dots,m$, let $A_i$ be an abelian variety with complex multiplication by $E_i$, and let $\sA_i$ be a class of varieties which dominates $A_i$.
Then $A = \prod_{i=1}^m A_i$ is dominated by $\prod_{i=1}^m \sA_i$.
\end{proposition}

\begin{proof}
Let $E$ be the compositum of the $E_i$.
Let $\sG = \Gal(\overline{E}/\bQ)$, and, $\sG_i = \Gal(\overline{E_i}/\bQ)$, where bars denote Galois closure.
Then, the linear disjointness of the $\overline{E_i}$ implies that $\sG = \prod_i \sG_i$.

From $\sG = \prod_i \sG_i$, and, the explicit description of the cocharacter groups of the Mumford-Tate groups of abelian varieties of CM-type given by Deligne \cite{900}*{p.~47}, we immediately infer that the special Mumford-Tate group of $A$ is the product of the special Mumford-Tate groups of $A_i$.
Thus $G = \prod_i G_i$, where $G$ is the Hodge group of $A$, and $G_i$ the Hodge group of $A_i$.

Let $V$ be an irreducible Hodge structure in the cohomology of $A$.
Let $W$ be an irreducible $G_{\bC}$-submodule of $V_{\bC}$.
Write $W = \bigotimes_i W_i$, with $W_i$ an irreducible $G_{i,\bC}$-module.
Since $\sG = \prod_i \sG_i$,
the set of Galois conjugates of $W$ consists of representations of the form
$W_1^{\sigma_1} \otimes \dots \otimes W_m^{\sigma_m}$,
with $\sigma_i \in \sG_i$.
Let $V_i$ be the direct sum of all $\sG_i$ conjugates of $W_i$.
Then, $V_{\bC}$ contains $V' = V_1 \otimes V_2 \otimes \dots \otimes V_m$.
Since $V'$ is $\sG$-invariant, it follows that $V'= V_0 \otimes E$ for some $G$-module $V_0$.
The irreducibility of $V$ then implies that $V = V_0$.

Now, each $W_i$ appears in $H^{a_i}(A_i,\bC)$ for some $a_i$.
Hence $W_i$ occurs in a Hodge structure $M_i \subset H^{b_i} (B_i, \bQ)$ for some $B_i \in \sA_i$, and such that $M_i$ contains $(b_i,0)$-forms.
Let $M = \bigotimes_i M_i \subset H^n(B,\bQ)$, where $n=\sum_i b_i$, and, $B=\prod_i B_i$.
Then, $M$ is equivalent to $V$ as a $G$-module, and contains $(n,0)$-forms.
\end{proof}

\section{Examples}

We give here all examples known to us of the general Hodge conjecture for CM abelian varieties.
This supersedes the listing in \cite{Abdulali2001}*{\S 4}.
We recall that a CM abelian variety is called \emph{nondegenerate} if its dimension equals the dimension of its Hodge group.
It is known (Hazama \cite{Hazama1983} and Murty \cite{Murty}) that if $A$ is a nondegenerate CM abelian variety, then, the Hodge ring of each power of $A$ is generated by divisors, so that the usual Hodge conjecture holds for all powers of $A$.

\begin{example}
\label{Tankeev}
Let $F$ be any totally real number field of degree $d$ over $\bQ$.
Then there exists a totally imaginary quadratic extension $E$ of $F$ such that $[\overline{E} : \overline{F}] = 2^d$, where bars denote Galois closure (Shimura \cite{ShimuraCanonical}*{1.10.1, p.~155}).
$\Gal(\overline{E}/\bQ)$ acts transitively on the set of CM-types of $E$ (Tankeev \cite{Tankeev2}*{Lemma~2.5, p.~187}).
Therefore, up to isogeny, there is only one abelian variety $A$ with complex multiplication by $E$.
Any power of $A$ is then dominated by the set of powers of $A$.
Since $A$ is nondegenerate (Ribet \cite{Ribet}*{Corollary 3.6, p.~87}), the general Hodge conjecture is true for all powers of $A$.
The general Hodge conjecture for $A$ was first proved by Tankeev \cite{Tankeev2}*{Theorem~2, p.~180}.
\end{example}

\begin{example}
If $F=\bQ$, and, $K$ is an imaginary quadratic number field then the hypotheses of the previous example are satisfied.
Therefore, the general Hodge conjecture is true for any power of an elliptic curve with complex multiplication.
This result is due to Shioda \cite{Shioda}.

Let $K_1$, $K_2$, $K_3$ be distinct quadratic imaginary number fields.
Then it is easy to see that they are linearly disjoint.
Let $E_i$ be an elliptic curve with complex multiplication by $K_i$ for $i=1,2,3$.
Then Prop.~\ref{product} implies the general Hodge conjecture for all $E_1^j \times E_2^k \times E_3^{\ell}$.
\end{example}

\begin{example}
We will now see that the general Hodge conjecture holds for any power of an abelian surface of CM-type.
Let $F$ be a totally real quadratic number field, and, $E$ a totally imaginary quadratic extension of $F$.
Let $A$ be an abelian variety with CM by $E$.
If $A$ is not simple, then it is a product of two elliptic curves, and as in the previous example, we have the general Hodge conjecture for all powers of $A$.

Suppose now that $A$ is simple.
If $E$ is Galois over $\bQ$, then $\sG := \Gal(E/\bQ)$ is cyclic of order $4$, acting transitively on the set of $4$ CM-types associated with $E$ (see \cite{ShimuraTaniyama}*{pp.~64--65}).
Thus, up to isogeny, $A$ is the only abelian variety with CM by $E$.
Therefore, $A$ is dominated by the set of powers of itself, and the general Hodge conjecture holds for all powers of $A$.

If $E$ is not Galois over $\bQ$, then, up to isogeny, there are two abelian surfaces with CM by $E$ (Moonen and Zarhin \cite{MoonenZarhin1999}*{p.~725}).
Call them $A$ and $B$.
Any power of $A$ is dominated by the set of all abelian varieties of the form $A^i \times B^j$.
But the Hodge group of $A \times B$ equals the product of the Hodge groups of $A$ and $B$ \cite{MoonenZarhin1999}*{Prop.~4.2, p.~725}.
Hence the usual Hodge conjecture holds for all $A^i \times B^j$.
Thus, the general Hodge conjecture is true for any power of an abelian surface of CM-type.
\end{example}

\begin{example}
Let $A$ be a simple abelian variety of CM-type $(E,\Phi)$ such that $d := \dim A$ is odd.
Let $F$ be the maximal real subfield of $E$, and $\overline{F}$ its Galois closure.
If $\Gal(\overline{F} / \bQ)$ is isomorphic to either the symmetric group or the alternating group on $d$ letters, then $A$ is nondegenerate (Dodson \cite{Dodson2}*{Prop.~2.1, p.~58}).
Suppose that $E$ does not contain any imaginary quadratic field.
Then, $[\overline{E}/\overline{F}] = 2^d$ (Dodson \cite{Dodson1}*{Prop.~2.2.2, p.~82}), so $A$ satisfies the hypotheses of Example~\ref{Tankeev}, and the general Hodge conjecture holds for all powers of~$A$.
\end{example}

\end{document}